\documentclass[10pt,letterpaper]{article}
\usepackage[top=0.85in,left=1.0in,footskip=0.75in,marginparwidth=2in]{geometry}
\usepackage{cite}
\usepackage{amsmath,amssymb,amsfonts}
\usepackage{amsthm}
\usepackage{algorithmic}
\usepackage{graphicx}
\usepackage{textcomp}
\usepackage{xcolor}
\usepackage{mathrsfs}
\usepackage{footmisc}

\newtheorem{proposition}{Proposition}
    
\usepackage[utf8]{inputenc}

\usepackage{cite}

\usepackage{nameref,hyperref}

\usepackage[right]{lineno}

\usepackage{microtype}
\DisableLigatures[f]{encoding = *, family = * }

\raggedright
\usepackage{changepage}

\usepackage[aboveskip=1pt,labelfont=bf,labelsep=period,singlelinecheck=off]{caption}

\makeatletter
\renewcommand{\@biblabel}[1]{\quad#1.}
\makeatother

\usepackage{color}

\definecolor{Gray}{gray}{.25}

\usepackage{graphicx}

\usepackage{sidecap}

\usepackage{wrapfig}
\usepackage[pscoord]{eso-pic}
\usepackage[fulladjust]{marginnote}
\reversemarginpar

\begin{document}
\vspace*{0.35in}

% title goes here:
\begin{flushleft}
{\Large
\textbf\newline{Adaptive Clutter Suppression via Convex Optimization}
}
\newline
% authors go here:
\\
Yifan He\textsuperscript{1},
Griffin Kearney\textsuperscript{2},
Makan Fardad\textsuperscript{1},
\\
\bigskip
\bf{1} Electrical Engineering \& Computer Science, Syracuse University, USA
\\[0.25cm]
\bf{2} Hidden Level Inc., Syracuse NY, USA
\\
\bigskip

\end{flushleft}

\section*{Abstract}
Passive and bistatic radar systems are often limited by strong clutter and direct-path interference that mask weak moving targets. 
Conventional cancellation methods such as the extensive cancellation algorithm require careful tuning and can distort the delay–Doppler response. 
This paper introduces a convex optimization framework that adaptively synthesizes per-cell delay–Doppler filters to suppress clutter while preserving the canonical cross-ambiguity function (CAF). 
The approach formulates a quadratic program that minimizes distortion of the CAF surface subject to linear clutter-suppression constraints, eliminating the need for a separate cancellation stage. 
Monte Carlo simulations using common communication waveforms demonstrate strong clutter suppression, accurate CFAR calibration, and major detection-rate gains over the classical CAF. 
The results highlight a scalable, CAF-faithful method for adaptive clutter mitigation in passive radar.

% the * after section prevents numbering
\section{Introduction}
Modern radar detection is frequently interference-limited: direct-path leakage, multipath, and stationary or slowly varying clutter can obscure weak moving targets, especially in passive and multistatic systems that exploit illuminators of opportunity (e.g. broadcast, cellular, satellite). 
Classical processing assumes availability of a trustworthy reference, performs disturbance suppression upstream, and then forms a delay-Doppler (DD) map by correlating the surveillance channel with time-frequency shifted replicas of the reference -- the cross-ambiguity function (CAF) \cite{Richards2014,GriffithsBaker2017}.\\

A large passive-radar literature tackles the ``cancel first, detect later'' paradigm via Extensive Cancellation Algorithms (ECA) and variants that project the surveillance data onto the orthogonal complement of a reference-derived interference subspace before CAF evaluation \cite{Yu2016,Palmarini2015,Colone2016,Zuo2021,Ansari2016}. 
These approaches can be effective, but they introduce tuning burdens (batch length, overlap, update rate) and are sensitive to reference-model mismatch; 
sliding or batch choices can Doppler-modulate slow targets or leave residuals if ill-tuned, motivating numerous refinements \cite{Palmarini2015,Colone2016,Zuo2021}.\\

In this paper we adopt a different strategy that avoids prescriptive upstream assumptions about the reference structure, providing an important advantage for passive radar where the illuminator’s modulation and bandwidth can evolve over time \cite{GriffithsBaker2017}. 
Instead of excising clutter and then correlating with simple time-frequency shifted reference signals, we synthesize, for each DD test cell, a linear filter whose coefficients are obtained by solving a convex quadratic program with linear constraints.
Crucially, the objective is not generic output-power minimization: it is to minimize distortion of the DD response surface relative to the canonical CAF surface, thereby preserving the interpretability and calibration of the DD map. 
Simultaneously, the linear constraints impose strong interference suppression on surveillance cells associated with clutter and direct-path/multipath, preventing target masking in each bin. 
In effect, the simple replica at a given DD is locally perturbed just enough to enforce clutter immunity while keeping the global surface as close as possible (in a convex least-squares sense) to what the unmodified CAF would have produced. 
The resulting problem is a convex quadratic program, solvable reliably with standard methods and readily batched across cells \cite{BoydVandenberghe2004}.\\

Conceptually, this is a constrained mismatched-filter design on the DD manifold that is aligned with the detection objective: we preserve the CAF’s desired response (in the optimization objective) while carving nulls where clutter resides (in the optimization constraints). 
The philosophy echoes well-studied ideas in adaptive/robust receiver design (e.g. MVDR/Capon beamforming and space-time adaptive processing a.k.a. STAP) which use linear constraints to maintain distortionless response while shaping nulls against interference \cite{Ward1994}; the novelty of our approach is to apply these ideas to the DD filter bank so that the CAF surface is preserved up to the minimum necessary distortion to ensure clutter immunity.\\

The proposed scheme (i) eliminates the need for a cancellation pre-processing stage, and (ii) produces a DD map that remains metrically comparable to the canonical CAF, facilitating CFAR calibration and operator interpretation. 
The formulation nests classical processing as a special case. 
With no clutter constraints, the optimizer reproduces the CAF; but with constraints active, it returns the closest CAF-like DD surface that also nulls designated clutter cells. 
The framework is readily extensible: regularization, inequality constraints for spectral or spatial notches, and multi-cell couplings for CFAR-like behavior fit naturally within the same convex program \cite{BoydVandenberghe2004,Richards2014} and are appealing areas for future work.\\

\paragraph*{Contributions}
\begin{enumerate}
    \item A CAF-faithful, constraint-driven convex optimization for per-cell adaptive filter synthesis. 
    The objective minimizes DD-surface distortion relative to the classical CAF, while equality constraints enforce strong clutter mitigation in surveillance cells.
    \item A quadratic program with linear equalities that enables efficient batched solutions and clean KKT interpretations.
    \item A performance characterization of the proposed method using Monte Carlo simulations involving random clutter and target scenarios.
\end{enumerate}

\section{Modeling the System}

Denote the discrete time-series reference waveform as $x$, a vector with $T$ components. Define $X$ as a $T \times N$ matrix with columns that are each a time/frequency shifted copy of the reference waveform $x$; here $N$ denotes the number of points in the DD surface that we wish to test. We note that in general $T \gg N$ so that $X$ is a very tall matrix.\\

We consider a data capture $y$ of the form
\begin{equation}
\label{eqn:y-expand}
y = X \rho + z,
\end{equation}
where $X$ is the matrix of time/frequency shifted copies of $x$, $\rho$ is the true DD response, and $z$ is interference modeled as a  vectorized white stochastic process with $\mathcal{E}(z) = 0$ and $\mathcal{E}(z z^*) = I$. Here $\mathcal{E}$ denotes the expectation operator and $z^*$ denotes the Hermitian (i.e., complex-conjugate transpose) of $z$.\\

At its core, DD processing is a method for estimating $\rho$.
In the classical approach, the estimate $\hat{\rho}$ is computed using 
\begin{equation}
    \label{eqn:dd-classic}
    \hat{\rho} = X^* y = X^* X \rho + X^* z,
\end{equation}
which is implemented via the CAF.
If $X^* X$ is diagonal then this technique provides a ``clean'' estimate up to random interference, but this is not the case in general.
The structure of $X^* X$ allows strong reflections from the environment (i.e., clutter) to mask the substantially weaker target reflections, which greatly degrades detection performance.\\

Clutter reflections correspond to a subset of the columns of $X$.
Without loss of generality we partition $X$ into surveillance and clutter columns, respectively $X_s$ and $X_c$, such that
\[
X =
\!\left[\begin{matrix}
X_s \!&\! X_c
\end{matrix}\right]\!,
\]
and this induces a parallel partition of the components of $\rho$,
\[
\rho =
\!\left[\begin{matrix}
\rho_s \\[0.1cm]
\rho_c
\end{matrix}\right]\!.
\]
In contemporary approaches, such as ECA, the data capture $y$ is preprocessed to mitigate clutter reflections before applying traditional DD estimation using (\ref{eqn:dd-classic}).
In this line of reasoning clutter cancellation is used to modify the data capture  yielding $\hat{y}$, a clutter-free signal of the form
\begin{equation}
    \label{eqn:y-eca}
    \hat{y} = y - X_c \rho_c = X_s \rho_s + z.
\end{equation}
In practice, the estimation of $X_c \rho_c$ can be difficult and may require signal-specific tuning to avoid target cancellation and contamination via interference.\\

We prescribe a new approach where, instead of modifying the data signal via cancellation pre-processing, we estimate $\hat{\rho}$ using an ``adapted'' matrix $U$ which is built through an optimization scheme.
The clutter signals (characterized by the columns of $X_c$) are prespecified, and linear equality constraints are used to firmly prevent target masking.
Independent of this, our optimization based approach affords us the freedom to prescribe objective functions that incentivize appealing characteristics.
In the present work we develop an objective that balances incentives derived from three distinct facets: reducing susceptibility to random interference, minimizing distortion of the nominal CAF, and improving auto-ambiguity artifacts present in the DD surface.
Pleasingly, we are able to motivate this treatment from first principles and intuitive arguments which ultimately yields a (convex) quadratic program with linear constraints.

\section{Delay-Doppler Processing via\\ Convex Optimization}

In an optimization approach it is tempting to find $\hat{\rho}$, an optimal estimate for $\rho$, by solving the least-squares problem of minimizing $\| y - X \rho \|^2$. This can be shown to be equivalent to computing 
\begin{align*}
\hat{\rho} 
&= U^* y
\\
&= U^* X \rho + U^* z,
\end{align*}
where $U$ is a matrix of smallest $\| U \|_\mathrm{F}^2$ that satisfies $U^* X = I$, and $\| \cdot \|_\mathrm{F}$ denotes the Frobenius norm. The rationale here is that $U^*$ is a left-inverse of $X$ and minimizes the total variance in $\hat{\rho}$
\begin{align*}
\mathrm{trace} \, \mathcal{E} \big( (U^* z)(U^* z)^* \big)
&= \mathrm{trace} \big( U^* \mathcal{E} (z z^*) U \big)
\\
&= \| U \|_\mathrm{F}^2.
\end{align*}

We next argue that the standard least-squares solution outlined above may not be the best approach to this problem. As done in the previous section, suppose we reorder the entries of $\rho$ such that it can be partitioned as
\[
\rho =
\!\left[\begin{matrix}
\rho_s \\[0.1cm]
\rho_c
\end{matrix}\right]\!,
\]
where $\rho_s$ represents the {\em surveillance} signal we seek (corresponding to reflections from nonnzero-doppler moving objects) and $\rho_c$ represents {\em clutter} interference (corresponding to reflections from zero-doppler static objects). As before, this reordering and partitioning in $\rho$ induces a similar structure in $X$ and $U$,
\[
X =
\!\left[\begin{matrix}
X_s \!&\! X_c
\end{matrix}\right]\!,
~~~
U =
\!\left[\begin{matrix}
U_s \!&\! U_c
\end{matrix}\right]\!,
\]
so that $\hat{\rho} = U^* X \rho + U^* z$ becomes
\begin{align*}
\left[\begin{matrix}
\hat{\rho}_s \\[0.1cm]
\hat{\rho}_c
\end{matrix}\right]\!
&=
\!\left[\begin{matrix}
U_s^* X_s \!& U_s^* X_c \\[0.1cm]
U_c^* X_s \!& U_c^* X_c
\end{matrix}\right]\!
\!\left[\begin{matrix}
\rho_s \\[0.1cm]
\rho_c
\end{matrix}\right]\!
+
\!\left[\begin{matrix}
U_s^* \\[0.1cm]
U_c^*
\end{matrix}\right]\! z.
\end{align*}
Using the fact that our interest is only in moving objects, we ignore the clutter block of the above equations and focus on the signal block
\begin{align}
\hat{\rho}_s = U_s^* X_s \rho_s + U_s^* X_c \rho_c + U_s^* z.
\label{eqn:rho-hat-s}
\end{align}
Ideally, we would like (a) $U_s^* X_s = I$ so that we recover $\rho_s$, (b) $U_s^* X_c = 0$ to prevent the corruption of the surveillance signal by clutter, and (c) the smallest value of 
$\| U_s \|_\mathrm{F}^2 = \mathrm{trace} \, \mathcal{E} \big( (U_s^* z)(U_s^* z)^* \big)$ to minimize the effect of interference as measured by its total variance. However, it turns out that the exact enforcement of $U_s^* X_s = I$ can result in a minimizer of $\| U_s \|_\mathrm{F}^2$ that performs poorly for DD processing.\footnote{To elaborate on this, and thinking of $U_s^* X_c$  as a matrix populated by inner products between the columns of $U_s$ with those of $X_s$, if the columns of $X_s$ are strongly aligned (or stated differently, they are roughly colinear) then satisfying $U_s^* X_s = I$ (which demands that the columns of $U_s$ be orthogonal to some columns of $X_s$ while having a unit inner product with others) can result in $U_s$ that has large columns and thus a large value of $\| U_s \|_\mathrm{F}^2$.}\\

We overcome the difficulty of large $\| U_s \|_\mathrm{F}^2$ by replacing the hard constraint $U_s^* X_s = I$ with a soft constraint in the context of an optimization problem; namely, we penalize in the objective function the deviation of $U_s^* X_s$ from the identity matrix.
However, this scheme can yield an extremely small $U_s$ if the weight of $\| U_s \|_\mathrm{F}^2$ is much larger than that of $\| X_s^* U_s - I \|_\mathrm{F}^2$ in the objective function.
We therefore include another soft penalty in the form of $\| U_s - X_s \|_\mathrm{F}^2$ to stabilize the scheme and to prevent excessive distortion of the DD surface.\\

We formulate the resulting optimization problem ($P$)
\[
\begin{array}{ll}
\text{minimize}
& 
\gamma_1 \, \| U_s \|_\mathrm{F}^2 + 
\gamma_2 \, \| U_s - X_s \|_\mathrm{F}^2 + 
\gamma_3 \, \| X_s^* U_s - I \|_\mathrm{F}^2
\\[0.15cm]
\text{subject to}
&
X_c^* U_s = 0
\end{array}
\]
where the values of $\gamma_i \geq 0$, $i=1,2,3$, $\gamma_1 + \gamma_2 + \gamma_3 = 1$, characterize the relative importance of the terms in the objective. This is a convex  program with a (convex) quadratic objective and linear constraints and thus admits a closed-form solution using the Lagrangian. It can be shown that, up to a scalar multiplicative factor, the solution of the above problem depends on $\gamma_1 + \gamma_2$ and not on the individual values of $\gamma_1$ and $\gamma_2$. In other words, for a fixed value of $\gamma_3$, the optimal $U_s$ changes only by a scalar multiple as $\gamma_1$ and $\gamma_2$ vary while their sum is held constant. This implies that the {\em relative} amplitude of the entries of $\hat{\rho}_s$, our estimated surveillance DD response in (\ref{eqn:rho-hat-s}), is only a function of $\gamma = \gamma_1 + \gamma_2$. Hence, moving forward, we set $\gamma_1 = 0$ for simplicity using the following Proposition as justification.\\

\begin{proposition}
\label{prop:1}
Let $\hat{U}_s(\gamma_1,\gamma_2,\gamma_3)$ be the optimal solution of ($P$) for $\gamma_1,\gamma_2,\gamma_3$. $\hat{U}_s$ is of the form
\[
    \hat{U}_s = 2 (\gamma_2 + \gamma_3) \, U_0(\gamma_1 + \gamma_2 , \gamma_3),
\]
where $U_0$ is a matrix with dimension equal to the dimension of $U_s$ and that depends only on $\gamma_1 +\gamma_2$ and $\gamma_3$.
\end{proposition}
We omit a detailed rigorous proof for brevity, instead summarizing the KKT conditions and the optimal solution.
Problem ($P$), as a convex quadratic program with linear equality constraints, is known to have a unique global solution that can be obtained by routine analytic methods through formulating the Lagrangian.
\begin{proof}
    Construction of the Lagrangian induces dual variables $\Lambda$, and we refer the reader to \cite{BoydVandenberghe2004} for a foundational treatment on the subject.
    The KKT conditions which constrain the optimal solution are given by
    \begin{align}
        \label{eqn:kkt_1}
        A \hat{U}_s + X_c \Lambda & = 2 (\gamma_2 + \gamma_3) X_s\\
        \label{eqn:kkt_2}
        X_c^* \hat{U}_s & = 0,
    \end{align}
    where
    \[
        A = 2\Big[ (\gamma_1 + \gamma_2) I + \gamma_3 X_s X_s^*].
    \]
    Notice that the matrix $A$ is a function of $\gamma_1 + \gamma_2$ and $\gamma_3$.
    The KKT conditions constitute a linear system which can be solved by standard methods and this yields
    \begin{equation}
        \label{eqn:optsolution}
        \hat{U}_s = 2 (\gamma_2 + \gamma_3) A^{-1} \Big( I - X_c (X_c^* A^{-1} X_c)^{-1} X_c^* A^{-1} \Big) X_s.
    \end{equation}
    The matrix
    \[
        A^{-1} \Big( I - X_c (X_c^* A^{-1} X_c)^{-1} X_c^* A^{-1} \Big)
    \]
    depends on the $\gamma$ parameters purely through $A$ and is therefore a function of $\gamma_1 + \gamma_2$ and $\gamma_3$.
    Thus take
    \[
        U_0 = A^{-1} \Big( I - X_c (X_c^* A^{-1} X_c)^{-1} X_c^* A^{-1} \Big) X_s,
    \]
    and the result is proved.
\end{proof}

Motivated by the preceding discussion, the main optimization problem used in this work to obtain $U_s$ is given by
\[
\begin{array}{ll}
\text{minimize}
& 
\gamma \, \| U_s - X_s \|_\mathrm{F}^2 + 
(1-\gamma) \, \| X_s^* U_s - I \|_\mathrm{F}^2
\\[0.15cm]
\text{subject to}
&
X_c^* U_s = 0
\end{array}
\]
where $\gamma \in [0,1]$. The resulting optimal $U_s$, together with $\hat{\rho} 
= U^* y$, renders our best estimate of $\rho_s$.\footnote{Note that $U_c$ can be set to zero or chosen arbitrarily, since we have no interest in $\rho_c$.}
We refer to the resulting DD surface as having been found through an {\em adaptive} scheme.

\section{Methods}
We simulate a communications waveform using standards defined in 3GPP TS 36.101 for the present study.
We used a standard reference measurement channel (RMC) waveform that was synthesized in MATLAB and which corresponds to the R.13 profile.
Passive radar built on these types of emissions are an active and appealing topic of current research \cite{Rai2021LTEPassiveRadarReview}.
This signal is used as the reference $x$ for testing and constitutes a representative $10$ MHz orthogonal frequency-division multiplexing (OFDM) transmission from a wireless base station to user equipment.
This waveform was selected as an interesting representative because of the prominence of OFDM signal architectures in the real-world.
The ubiquity of these signals supports their utility since infrastructure is widely available, especially in urban environments, to enable mobile electronics (e.g. smartphones).\\

We use the RMC signal to generate simulated data which contains randomly generated clutter, targets, and thermal noise.
We then process the simulated data using the adapted optimization scheme and the ``unadapted'' baseline in which the classical approach of (\ref{eqn:dd-classic}) is used for comparison.
The DD surface that is produced by the adapted and baseline data paths are post processed using a standard cell-averaged (CA) constant false-alarm rate (CFAR) detector \cite{Richards2014}.
And we ultimately evaluate the scheme by aggregating statistical performance over these batched random simulations.

\subsection{Simulation Design}
Simulated data is generated in $10$ ms windows, corresponding to $T = 153600$ with $f_s = 15.36$ MSPS.
We inject random complex-valued Gaussian data to model thermal noise in addition to modeled clutter and target returns in each simulated data capture.\\

\subsubsection{Clutter}
Discrete clutter is simulated by specifying an integer $n_c$ - denoting a number of discrete static clutter sources - which are associated with random radar cross sections (RCS) and time delays.
The RCS and time delay values are used to represent random clutter object sizes and distances from the simulated receiver, respectively.
Different random distributions for these values correspond to different types of physical scenarios.
We prescribe $n_c = 1000$ and a Gaussian distribution for the clutter RCS with mean $-6$ dB and standard deviation $6$ dB.
We prescribe clutter delays drawn from a uniform distribution on the interval $[0,\tau_c]$ for a specified max delay, $\tau_c$.
We set $\tau_c = 4.0365 \mu s$ which corresponds to a clutter range approximately between $0$ km and $1.2$ km.
The individual clutter responses are additively aggregated into a total clutter response and used in the performance evaluation.\\

This prescribed clutter model is simple in nature and serves as a baseline to demonstrate preliminary validation of the new optimization scheme.
State-of-the-art clutter simulations could be implemented in future work, such as the richer models described in \cite{Gogineni2022}, but this simple model is sufficient for initial baseline evaluation of the proposed method.\\

\subsubsection{Target}
Random target returns are generated in simulated data by applying delay and frequency shifts to the reference signal which are drawn from a uniform probability distribution.
Targets were simulated with Doppler shift magnitudes (both positive and negative) between $300$Hz and $1300$Hz and time delays between $0.39 \mu s$ and $3.71 \mu s$.
The target return is scaled according to a randomly selected target SNR modeled with a Gaussian distribution with mean $2$ dB and standard deviation $0.33$ dB.

\subsection{Target Detection \& False-Alarms}
The DD surface is a lattice of test responses, represented as an $N_{doppler} \times N_{delay}$ matrix we denote as $S$, which is constructed by reshaping the estimated $\hat{\rho}_s$ vector.
Each row of $S$ corresponds to the response of the observation data to a fixed frequency (Doppler) shift of the reference across time delays.
Similarly, each column of $S$ corresponds to the response of the observation data to a fixed time delay of the reference across Doppler shifts.\\

The CFAR detector is used to process $S$ and it produces a boolean result for each test cell with $1$ indicating that a target is present and $0$ indicating that a target is absent.
The matrix of underlying delay and Doppler shifts that are used to generate $S$ discretize the domain.
However, in simulation (as in the real-world) we do not require that radar targets are endowed with signal responses that fall within the discrete set of shifts in general.
This leads to the excitation of neighboring cells in $S$ corresponding to the nearest bounding delay and Doppler shifts.
So long as the CFAR detector accurately flags at least one of the four neighboring target cells, then we treat the target as correctly detected.\\

\begin{figure}
    \centering
    \includegraphics[width=\linewidth]{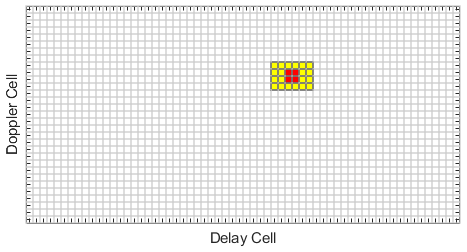}
    \caption{Detection and false-alarm scoring schematic. An example target neighborhood (red) and exclusion neighborhood (yellow) is depicted for illustrative purposes.}
    \label{fig:1}
\end{figure}

In parallel with target detection we also track false-alarms, when the CFAR detector indicates $1$ at DD cells that are not near a true target.
Waveform artifacts cause additional spread in target responses which complicates this matter.
For instance, the RMC waveform we use for investigative purposes induces structured smearing in the delay dimension.
To accommodate waveform artifacts without unnecessarily penalizing the detection scheme we include exclusion neighborhoods that superset the target neighborhoods.
Any positive CFAR test cells in the exclusion neighborhood are not counted when computing the false-alarm rate.
We illustrate the scheme in Figure \ref{fig:1}.
We evaluate false-alarm performance using an exclusion neighborhood which is $4$ test cells by $6$ delay cells.
The CFAR detector was configured using $1$ guard cell in both delay and Doppler dimensions, with $12$ training cells in delay and $6$ training cells in Doppler around each cell-under-test (CUT).
We configure the desired probability of false-alarm across different values to survey the performance of the scheme under various CFAR tunings.

\section{Results \& Illustrative Example}
We simulated $10000$ signal measurements as described above and used it for exploration of various adaptive tuning schemes.
We evaluated the adaptive processing by examining detection and observed false-alarm rates as a function of the tuned CFAR false-alarm rate and the parameter $\gamma$.
Throughout all simulated test cases we define the clutter space $X_c$ as the set of zero-Doppler, \emph{whole sample} (as opposed to fractional) time delays from $0$ to $4.0365 \mu s$.\\

\subsection{Illustrative Performance}
\begin{figure}
    \centering
    \includegraphics[width=\linewidth]{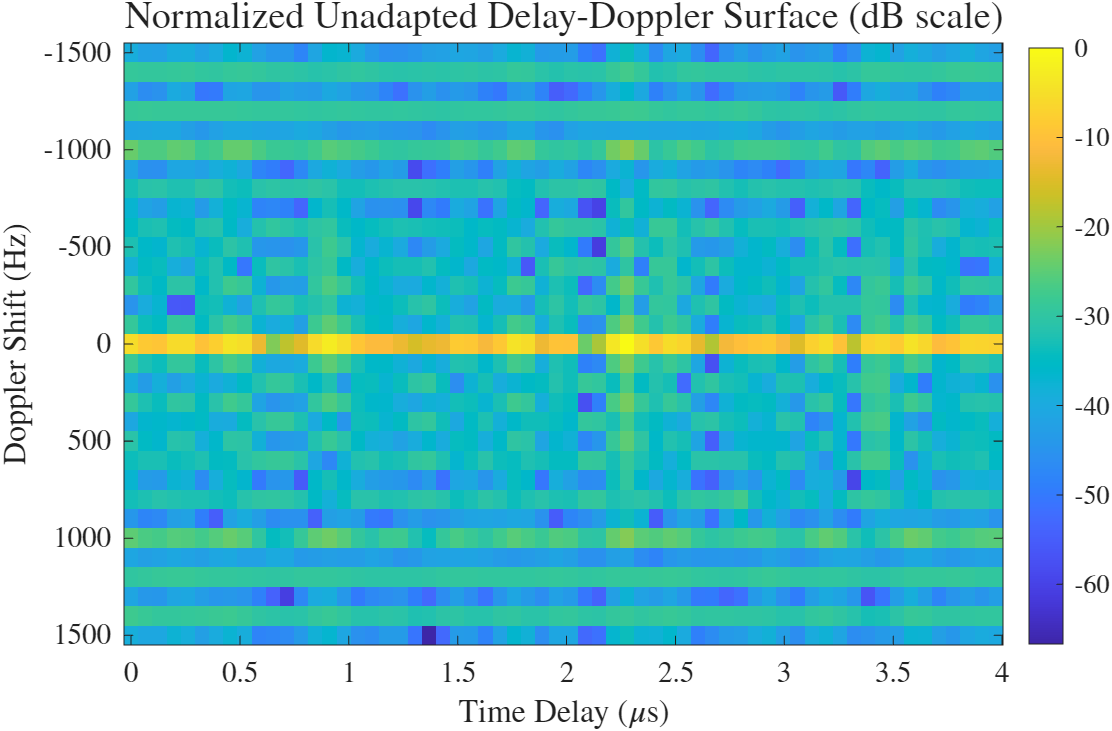}
    \caption{Example normalized unadapted DD plot.}
    \label{fig:2}
    \vspace{0.3cm}
    \includegraphics[width=\linewidth]{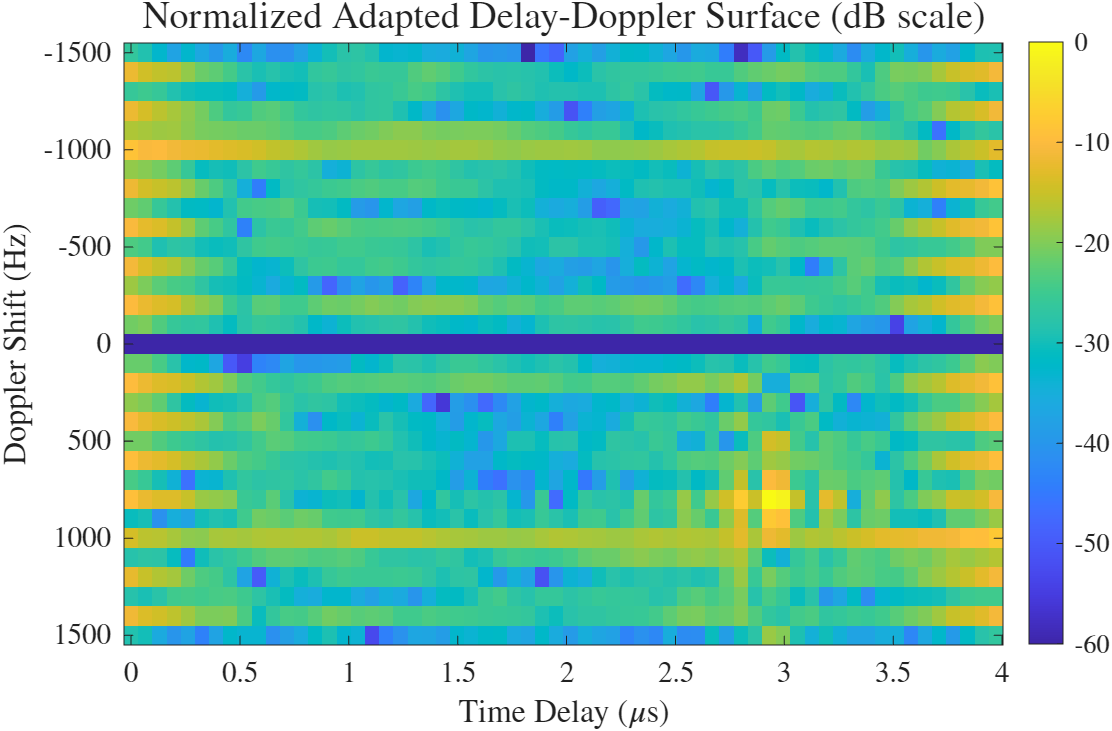}
    \caption{Example normalized adapted DD plot.}
    \label{fig:3}
\end{figure}

First we show a comparison of unadapted and adapted DD surfaces using a representative simulated data capture in Figures \ref{fig:2} and \ref{fig:3}, respectively.
In each, the DD surface is normalized by the maximal magnitude across the test cells which scales the highest peak of the surface to $0$ dB for comparability.
It is clear that the DD response is dominated by clutter in the unadapted baseline.
However, after applying the adaptive processing the clutter is heavily suppressed which greatly improves target visibility.
This comparison demonstrates the effectiveness of the adaptive scheme in eliminating clutter masking of targets. 
Intrinsic to the scheme, the adaptive parameter $\gamma$ must be set before processing can occur.
We show in the next section that tuning this parameter to values slightly less than one -- in other words prescribing a slight incentive to enforce surveillance filter orthogonality -- leads to the most appealing performance.\\

\subsection{Adaptive Parameter Tuning}
\begin{figure}[h!]
    \centering
    \includegraphics[width=\linewidth]{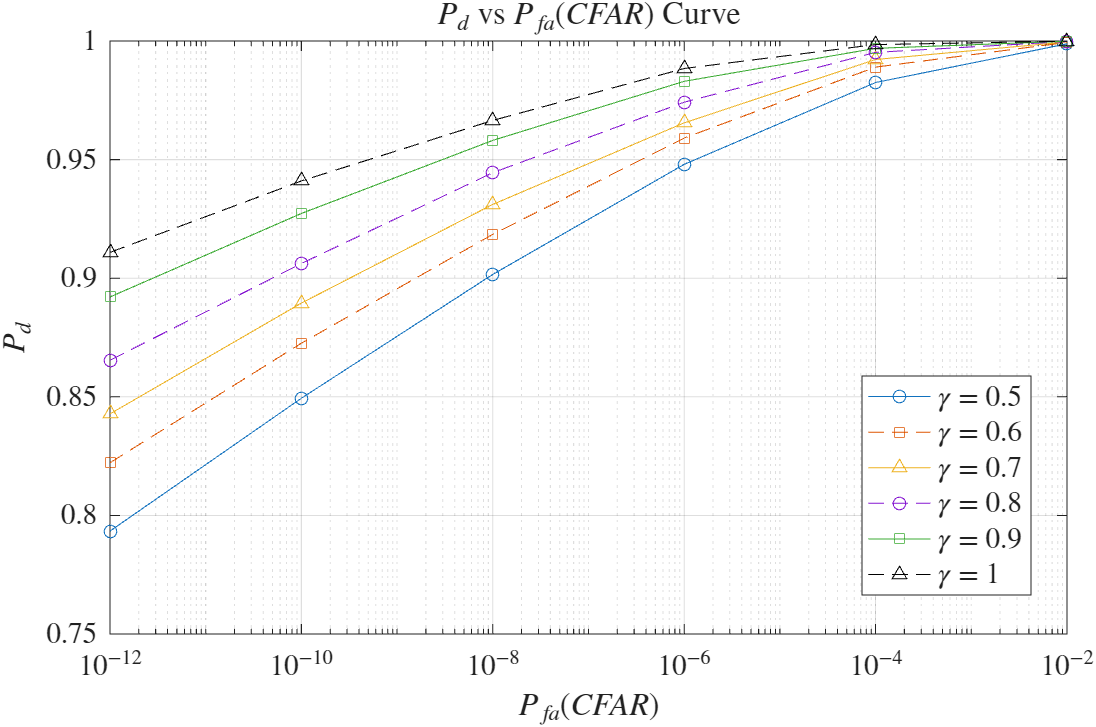}
    \caption{Coarse $P_d$ vs $P_{fa} (CFAR)$ $\gamma$ tuning curves.}
    \label{fig:4}
\end{figure}

\begin{figure}[h!]
    \centering
    \includegraphics[width=\linewidth]{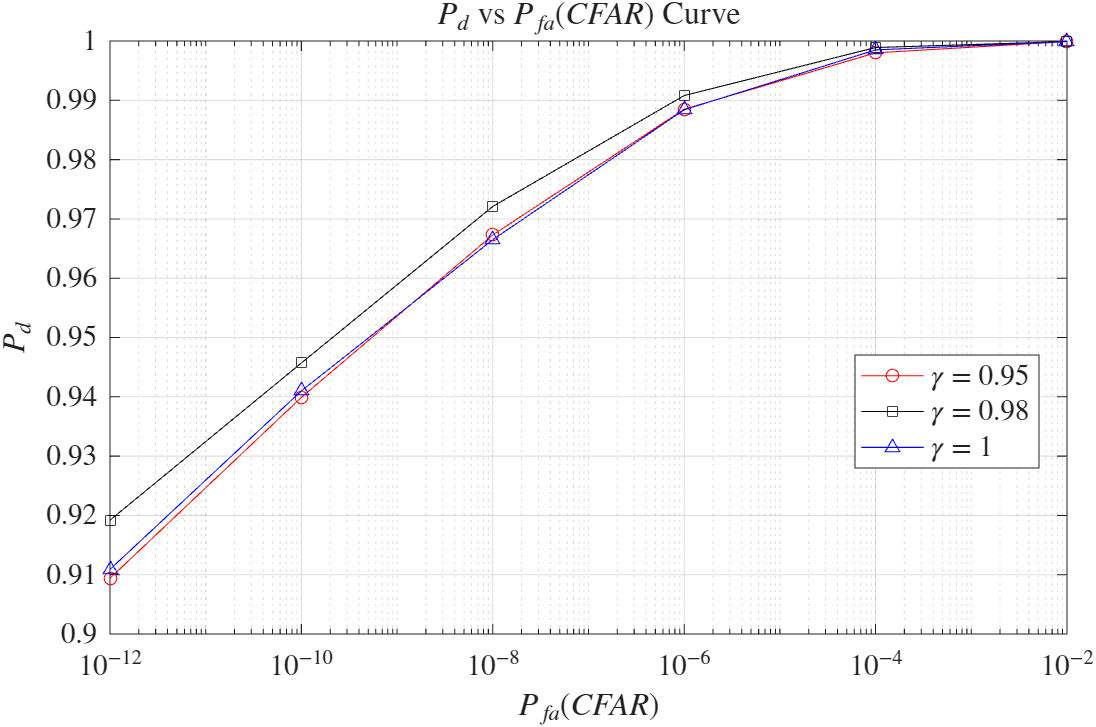}
    \caption{Fine $P_d$ vs $P_{fa} (CFAR)$ $\gamma$ tuning curves.}
    \label{fig:5}
\end{figure}

We explore both a coarse and a fine tuning of the adaptive parameter $\gamma$.
We evaluate each configuration by examining the receiver operating characteristic (ROC) curve where the target detection rate is plotted against the \emph{tuned} CFAR false-alarm rate.
We prescribe a coarse set of adaptive test parameters in Figure \ref{fig:4} which shows that the performance of the scheme decays with substantial reductions in $\gamma$.
In contrast, Figure \ref{fig:5} shows an examination of a fine set of adaptive test parameters near $\gamma = 1$ and demonstrates that performance is benefited by slight deviations from unity.
Fortunately, the results are not notably sensitive near $\gamma = 1$. 
We consistently observe substantial improvement and thus in the absence of exhaustive testing and calibration the adaptive parameter can be set naively to unity with only minor performance sacrificed.
We found that $\gamma = 0.98$ exhibited good performance in the test set we investigated.
We compare both the detection and observed false-alarm rates of this optimal tuning against the unadapted baseline in Table \ref{tab:1}.\\

In summary, if $\gamma$ is set too low, then it causes performance degrading distortion in the DD surface, but slight reductions in $\gamma$ (i.e., values of $\gamma$ less than but close to unity) provide reduced false-alarm rates by mitigating auto-ambiguity artifacts of the R.13 reference signal.
In practice, the optimal $\gamma$ value may differ depending on reference signal characteristics and other specifics of the particular case at hand.
\begin{table}[h!]
    \centering
    \begin{tabular}{|c|c|c|c|c|}
        \hline
        \textcolor{white}{--} & \multicolumn{2}{c|}{$\mathrm{Unadapted}$} & \multicolumn{2}{c|}{$\gamma$ = $0.98$}\\
        \hline
         $P_{fa}$ (CFAR) & $P_d$ & $P_{fa}$ & $P_d$ & $P_{fa}$\\
        \hline
        $10^{-2}$  & 0.1627 &  0.0413 & 0.9999 & 5.42 $\times \mathrm{10}^{\mathrm{-2}}$\\
        \hline
        $10^{-4}$  & 0.0194 &  0.0208 & 0.9989 & 6.96 $\times \mathrm{10}^{\mathrm{-3}}$\\
        \hline
        $10^{-6}$  & 0.0033 &  0.0141 & 0.9908 & 5.46 $\times \mathrm{10}^{\mathrm{-6}}$\\
        \hline
        $10^{-8}$  & 0.0004 &  0.0100 & 0.9721 & 0 \\
        \hline
        $10^{-10}$  & 0.0000 & 0.0073 & 0.9457 & 0 \\
        \hline
        $10^{-12}$  & 0.0000 & 0.0053 & 0.9192 & 0 \\
        \hline
    \end{tabular}
    \vspace{0.05cm}
    \caption{Detection and false-alarm rate comparisons in Monte Carlo testing: unadapted vs adapted with $\gamma = 0.98$}
    \label{tab:1}
\end{table}

\section{Conclusions \& Future Work}
We have provided preliminary evidence that the proposed approach exhibits significant promise as an adaptive framework for radar processing.
Using a well-tuned adaptive parameter ($\gamma = 0.98$), we observed detection rates exceeding $90\%$ with well-controlled false-alarm rates under CFAR detection in simulations; in the unadapted baseline, detections were virtually nonexistent.
However, the present results are contingent on the accuracy of our simulation scheme as a representation of real data systems.
Future work is focused on two major thrusts (i) theoretical study and advancement of the optimization framework (ii) applicability and performance using real-world data.\\

In (i) it is important and beneficial to quantify how susceptible the new processing scheme is to reference errors (i.e. when the modeled $x$ differs from the truth) as they are prominent in passive radar processing systems where the illuminator of opportunity is uncooperative.
The robust optimization approach we described enables extensions to help protect the scheme from performance degradations due to errors, such as those through the inclusion of augmented constraints or modified objectives.\\

In parallel, the existing framework should be further validated through real RF data analysis in pursuit of (ii).
Syracuse, NY is an area with a rich history of RF and radar system development.
Joint projects between Syracuse-based leaders in the RF sensor industry (e.g. Hidden Level, Inc) and Syracuse University provide an avenue to achieve this goal.
The consideration of real-world systems will not only emphasize applicability of future development to the generation of real technologies, but it will lay bare the challenges that must be overcome to effectively deploy research techniques. 
This includes the study of real-world RF phenomenology that degrade system performance (e.g. reference contamination, extraneous interference) across diverse sets of passive illuminators and efficient implementations to enable necessary real-time processing.

%This is where your bibliography is generated. Make sure that your .bib file is actually called library.bib
\bibliography{ieee_radar_conf_2026.bib}

%This defines the bibliographies style. Search online for a list of available styles.
\bibliographystyle{abbrv}

\end{document}